\documentclass[reqno]{amsart}
\usepackage{tikz}
\usetikzlibrary{patterns} 
\usepackage{pgfplots}
\usepackage{amsfonts,amsmath,amsthm,amssymb,mathrsfs,ctex}
\usepackage{color}
\usepackage{amsmath,amssymb,amsfonts, bm}
\usepackage{graphicx}
\usepackage{newunicodechar}
\usepackage{xcolor}
\usepackage{geometry}
\usepackage{enumitem}
\numberwithin{equation}{section}
\usepackage[pdfstartview=FitH,colorlinks=true]{hyperref}
\hypersetup{urlcolor=red, linkcolor=blue,citecolor=red}
\allowdisplaybreaks

\theoremstyle{plain}

\usepackage{ctex}
\usepackage{bm}
\usepackage{bbm}
\newtheorem{theorem}{Theorem}[section]
\newtheorem{definition}{Difinition}[section]
\newtheorem{lemma}[theorem]{Lemma}
\newtheorem{proposition}{Proposition}[section] 
\newtheorem{remark}[theorem]{Remark}

\newlist{romanlist}{enumerate}{1} 
\setlist[romanlist]{label=(\roman*), leftmargin=*}
\newlist{theoremlist}{enumerate}{1}
\setlist[theoremlist]{label=(\roman*), leftmargin=*}

\begin{document}
\title{Global Well-Posedness and Scattering for the Fractional Inhomogeneous Nonlinear Schrödinger Equation with Time-Dependent Damping}
\author[]
{Boyu Jiang, Kexue Li*}

\address{Boyu Jiang
 \newline\indent
 School of Mathematics and Statistics, Xi'an Jiaotong University, Xi'an ,710049, Shaanxi,China
 \newline\indent
 Kexue Li
 \newline\indent
 School of Mathematics and Statistics, Xi'an Jiaotong University,Xi'an 710049, Shaanxi,China
 \newline\indent
 }
 \email{ jiangboyu@stu.xjtu.edu.cn; kxli@mail.xjtu.edu.cn}
\date{}
\begin{abstract}
In this paper，we prove the global well posedness and scattering of solutions for a nonlinear Schrodinger type equation with fractional Laplacian, inhomogeneous nonlinearity and time-dependent damping.We give an elementary proof based on Strichartz estimates combined with estimates of the nonlinear term in mixed-norm spaces. To avoid the loss of regularity in the Strichartz estimates, we consider the radial data.
\end{abstract}
\maketitle
\section{Introduction}
In this paper, we mainly study the following fractional inhomogeneous nonlinear Schr\"odinger equation with time-dependent damping:
\begin{equation}\label{一般情形}
    \left\{
\begin{array}{ll}
\mathrm{i}\partial_t u - (-\Delta)^s u + ia(t)u = \mu |x|^{-b}|u|^{2\sigma}u, & \quad (t,x) \in (0,\infty) \times \mathbb{R}^N, \\[1ex]
u(0,x) = u_0(x)\in H^s(\mathbb{R}^N), & \quad x \in \mathbb{R}^N,
\end{array}
\right. 
\end{equation}
here \(\mu=\pm 1\), \(p>1\), and \(a\in C([0,\infty))\) is a nonnegative continuous function representing dissipation.

By the simple change of variables \(v(t,x)=\mathrm{e}^{A(t)}u(t,x)\), we can transform the equation into the following form:
\begin{equation}\label{变量替换之后的方程}
\begin{cases}
\mathrm{i}\partial_t v - (-\Delta)^s v = \mu\mathrm{e}^{(-2\sigma)A(t)}|x|^{-b}|v|^{2\sigma}v, \\
v(0, x) = u_0(x)\in H^s(\mathbb{R}^N).
\end{cases}
\end{equation}

Equation \eqref{一般情形} has profound applications in many physical contexts. For example, time-dependent damping is very useful in dealing with Bose--Einstein condensate (BEC) models; or in quantum mechanics, when the viscosity is a function of time, time-dependent friction also arises.

When the damping term \(a(t)\equiv0\), \(s=1\), and \(b=0\), equation \eqref{一般情形} reduces to the standard nonlinear Schr\"odinger equation. The study of well-posedness and scattering theory for the nonlinear Schr\"odinger equation is nearly complete; see mainly \cite{cazenave2003,Tao2006} and the references therein. When \(a(t)\) is a constant, many researchers have also studied this aspect, and global existence has been obtained in the defocusing subcritical case, i.e., \(\mu=1,\ 1<p<1+\frac{4}{(N-2)^+}\), regardless of the size of the damping; see \cite{cazenave2003}. However, in the focusing case, i.e., \(\mu=-1\), a global existence result is available only for sufficiently large damping and when the equation is energy-subcritical and mass-supercritical \((1+4/N\le p\le1+4/(N-2)^+)\); see \cite{ohta2009remarks}. Also in the focusing case, for mass-subcritical \(p<1+4/N\), the standard nonlinear Schr\"odinger equation admits blow-up solutions; when the size of the damping is relatively small, it still cannot prevent the blow-up effect induced by the nonlinear term, and blow-up solutions still exist. Only when the size of the damping term is sufficiently large can global existence be obtained; see \cite{hamouda2023existence}. It is worth mentioning that the method here can also be extended to the variable-coefficient damping case.

When \(b=0\) and \(a(t)\equiv0\), the equation reduces to the fractional Schr\"odinger equation. Considerable research has also been devoted to the fractional Schr\"odinger equation by many experts and scholars. Its earliest origin can be traced to \cite{Laskin2002Fractional}; later, Sire et al. studied local well-posedness and the corresponding Strichartz estimates in \cite{HongSire2015Fractional}. Subsequently, some experts and scholars studied global well-posedness and scattering on this basis; see, for example, \cite{Sire2023SingularLevy,Sun2018Scattering,Guo2018EnergyCritical}. When \(b\neq 0\), the equation is called the inhomogeneous nonlinear fractional Schr\"odinger equation, and many researchers have also made considerable studies in this direction; see, e.g., \cite{Radhia2024,Saanouni2022Local}.

We also give the following overview of the long-time behavior of the equation. When the equation lies in the subcritical regime, i.e., \(1<p<1+4/(N-2)^+\), for constant-coefficient damping, the scattering results for \eqref{一般情形} can be found in \cite{dinh2021blowup,inui2019asymptotic}; whether the equation is defocusing or focusing, provided suitable conditions are satisfied, the above authors all obtained satisfactory scattering results. Later, also in the subcritical case, for variable-coefficient damping, several authors studied scattering in a certain sense in \cite{hamouda2023existence}, assuming the condition \(\underline{\mathbf{a}}>0\) (see the definition in the next section). The authors proved a scattering result for the energy-subcritical and mass-supercritical equation in the focusing case, naturally building on the global existence result in \cite{hamouda2023existence}.

It is worth noting that, in view of the class of initial data we consider, the above assumption on the damping term is natural to some extent. In fact, we believe that there is a balance (constraint) relationship between the condition on the damping and the initial data. In this context, the damping \(t\mathbf{a}(t)\sim\gamma>0\) (as \(t\to\infty\)) considered in [3] requires a class of oscillatory initial data.

In this paper, we focus on the global existence and scattering results for the Schr\"odinger equation with variable damping under the energy-subcritical condition, in both the focusing and defocusing cases. Before beginning the main presentation of the paper, we first state the main theorems.
\begin{theorem}\label{main result}
   Assume that the initial value of equation \eqref{一般情形} belongs to \(H^s_{rd}\), where this denotes the radial Sobolev space. If \(\mu=-1\), under the condition \(0<\sigma<\frac{2s-b}{N-2s}\) (i.e., the subcritical case), provided that \(I(u_0)>0\) and \(E(u_0)<E(Q)\), the solution \(u(t)\) of \eqref{一般情形} is a global solution. When \(\mu=1\), for any radial initial value \(u_0\in H^s_{rd}\), the solution exists globally. On the other hand, assuming that a global solution exists and $\underline{a}>0$, then the solution scatters regardless of whether the equation is focusing or defocusing.
\end{theorem}

The paper is organized as follows. In Section 2, we give some main definitions and notation, and in Section 3 we prove our main theorems. Unless otherwise stated, in the sequel we use \(A\lesssim B\) to mean that there exists a constant \(C\) such that \(A\le CB\), and \(A\approx B\) to mean that \(A\lesssim B\) and \(B\lesssim A\).
\section{Some Main Definitions and Notation}
In this section, we are devoted to giving some important estimates and some necessary definitions. First, we introduce the following notation:
\begin{equation}\label{A(t)定义}
    A(t) = \int_{0}^{t} a(s) \mathrm{d}s, 
\end{equation}
\begin{equation}
    \underline{a} = \inf_{t>0} \left( \frac{A(t)}{t} \right). 
\end{equation}
Obviously, it can be seen that \(\underline{a}>0\) implies
\begin{equation}
    \int_{0}^{\infty} a(s) ds = \infty. 
\end{equation}
In addition, we write
\begin{equation*}
    \| \cdot \|_{L^r} := \| \cdot \|_{L^r(\mathbb{R}^N)}, \quad \| \cdot \| := \| \cdot \|_{L^2}, \quad \| \cdot \|_{H^1} := \| \cdot \|_{H^1(\mathbb{R}^N)}.
\end{equation*}
\begin{definition}
Let \(p\in[1,\infty]\). By \(L^p(I,X)\) we denote the set of (equivalence classes of) measurable functions \(f:I\to X\) such that the function \(t\mapsto\|f(t)\|\) belongs to \(L^p(I)\). For \(f\in L^p(I,X)\), its norm is defined by
\[
\|f\|_{L^p(I,X)} =
\begin{cases}
\displaystyle \left( \int_I \|f(t)\|^p \, dt \right)^{1/p}, & \text{if } p<\infty, \\[15pt]
\displaystyle \operatorname*{ess\,sup}_{t\in I} \|f(t)\|, & \text{if } p=\infty.
\end{cases}
\]
\end{definition}
Here \(X\) is often taken to be an \(L^r\) space.
\begin{definition}[Sobolev space]
Let \(m\in\mathbb{N}\) and \(1\le p\le\infty\). The Sobolev space \(W^{m,p}(\Omega)\) is defined by
\[
W^{m,p}(\Omega)
=
\{u\in L^p(\Omega): D^\alpha u\in L^p(\Omega)\ \text{for } |\alpha|\le m\}.
\]
Equipped with the norm \(\|\cdot\|_{W^{m,p}}=\|\cdot\|_{W^{m,p}(\Omega)}\) defined as follows, \(W^{m,p}(\Omega)\) is a Banach space:
\[
\|u\|_{W^{m,p}}
=
\sum_{0\le|\alpha|\le m}\|D^\alpha u\|_{L^p(\Omega)}.
\]

If \(p<\infty\), the closed subspace \(W_0^{m,p}(\Omega)\) of \(W^{m,p}(\Omega)\) is defined as the closure of \(\mathcal{D}(\Omega)\) in \(W^{m,p}(\Omega)\).

When \(p=2\), we write \(W^{m,2}(\Omega)=H^m(\Omega)\) and \(W_0^{m,2}(\Omega)=H_0^m(\Omega)\), and equip \(H^m(\Omega)\) with the equivalent norm
\[
\|u\|_{H^m(\Omega)}
=
\|u\|_{H^m}
=
\left(
\sum_{0\le|\alpha|\le m}\int_{\Omega}|D^\alpha u(x)|^2\,dx
\right)^{1/2}.
\]
Then \(H^m(\Omega)\) (and hence \(H_0^m(\Omega)\)) is a Hilbert space with inner product
\[
(u,v)_{H^m}
=
\sum_{0\le|\alpha|\le m}
\operatorname{Re}\int_{\Omega}
D^\alpha u(x)\overline{D^\alpha v(x)}\,dx.
\]
\end{definition}

For the sake of convenience in what follows, we also present the Sobolev embedding theorem in its general form, 
\begin{theorem}[Sobolev embedding theorem]
Suppose that \(\Omega\) has a Lipschitz continuous boundary. Then the following properties hold:
% using the newly defined romanlist environment
\begin{romanlist}
    \item If \(1 \le p < N\), then for any \(q \in \left[p, \frac{Np}{N-p}\right]\), we have \(W^{1,p}(\Omega) \hookrightarrow L^q(\Omega)\).
    \item If \(p = N > 1\), then for any \(q \in [p, \infty)\), we have \(W^{1,p}(\Omega) \hookrightarrow L^q(\Omega)\).
    \item If \(p = N = 1\), then for any \(q \in [p, \infty]\), we have \(W^{1,p}(\Omega) \hookrightarrow L^q(\Omega)\).
    \item If \(p > N\), then \(W^{1,p}(\Omega) \hookrightarrow L^\infty(\Omega)\).
\end{romanlist}

If \(\Omega\) has a uniformly Lipschitz continuous boundary, then:
\begin{romanlist}[start=5] 
    \item If \(p > N\), then \(W^{1,p}(\Omega) \hookrightarrow C^{0,\alpha}(\overline{\Omega})\), where \(\alpha = \frac{p-N}{p}\).
\end{romanlist}
\end{theorem}

Below we present the crucial Strichartz estimates for the fractional operator:
\begin{definition}[admissible pairs]
\begin{enumerate}
    \item If \(q,r\ge 2\) and satisfy the following conditions, then the pair \((q,r)\) is called admissible:
    \[
    \frac{4N+2}{2N-1} \leq q \leq \infty, \quad \frac{2}{q} + \frac{2N-1}{r} \leq N - \frac{1}{2},
    \]
    or
    \[
    2 \leq q \leq \frac{4N+2}{2N-1}, \quad \frac{2}{q} + \frac{2N-1}{r} < N - \frac{1}{2}.
    \]
    
    \item We denote by \(\Lambda_\alpha\) the set of admissible pairs \((q,r)\) such that \((N,q,r)\neq(2,2,\infty)\) and
    \[
    N\left(\frac{1}{2} - \frac{1}{r}\right) = \frac{2s}{q} + \alpha.
    \]
    
    \item If \(I\subset\mathbb{R}\) is an interval, we use the following notation for the Strichartz spaces and norms:
    \begin{equation}
    \Omega_\alpha(I) := \bigcap_{(q,r)\in\Lambda_\alpha} L^q(I,L^r), \quad \|\cdot\|_{\Omega_\alpha(I)} := \sup_{(q,r)\in\Lambda_\alpha} \|\cdot\|_{L^q(I,L^r)}. \tag{3}
    \end{equation}
    For brevity, we write \(\Lambda_0:=\Lambda\), \(\Omega(I):=\Omega_0(I)\), and \(\Omega_\alpha:=\Omega_\alpha(\mathbb{R}_+)\).
\end{enumerate}
\end{definition}
\begin{proposition}[Strichartz estimates]
Let \(N\ge 2\), \(\frac{2}{2N-1}<s\le 1\), and let \(f\) and \(g(t,\cdot)\) be radial functions. Then
\begin{enumerate}
    \item \(\|e^{-i(-\Delta)^s} f\|_{\Omega_\alpha(I)} \leq C \|(-\Delta)^{\frac{\alpha}{2}} f\|;\)
    \item \(\left\| \int_0^t e^{-i(t-t')(-\Delta)^s} g(t', \cdot) \, dt' \right\|_{\Omega_\alpha(I)} \leq C \|g\|_{\Omega_{-\alpha}(I)}.\)
\end{enumerate}
\end{proposition}
For the damped equation, we define the following three quantities (the first two of which are referred to as the mass and the energy, respectively).
\begin{align*}
\mathbf{M}(u(t)) &= \int_{\mathbb{R}^N} |u(t, x)|^2 \, \mathrm{d}x, \\
\mathbf{E}(u(t)) &= \|\nabla^s u(t)\|_{L^2}^2 + \frac{2\mu}{p+1} \int_{\mathbb{R}^N} |x|^{-b}|u(t, x)|^{p+1} \, \mathrm{d}x, \\
\mathbf{I}(u(t)) &= \|\nabla^s u(t)\|_{L^2}^2 +\mu \int_{\mathbb{R}^N} |x|^{-b}|u(t, x)|^{p+1} \, \mathrm{d}x.
\end{align*}

We have the following proposition, taken from \cite{hamouda2025energy}.
\begin{proposition}
Let \(u\) be a sufficiently smooth solution of equation (1.1) on \(0\le t\le T\). Then the following hold:
\begin{align}
M(u(t)) &= \mathrm{e}^{-2A(t)} M(u_0),  \label{质量耗散},\\
\frac{d}{dt}E(u(t)) &= -2a(t)I(u(t)).  \label{能量导数}
\end{align}
\end{proposition}
In what follows, we present the definition of scattering adopted in this paper, 
\begin{definition}[forward scattering]
We say that a global solution \(u\) of equation \eqref{变量替换之后的方程} scatters forward in \(H^1(\mathbb{R}^N)\) if there exists \(u^+\in H^1(\mathbb{R}^N)\) such that
\[
\lim_{t\to\infty}\left\| e^{it(-\Delta)^s}u(t)-u^+ \right\|_{H^s}=0.
\]
\end{definition}
\begin{remark}
   Note that the energy of the damped equation has a decay property; hence it is meaningless to discuss backward scattering for the damped equation. Therefore, we only give the definition of forward scattering.  
\end{remark}
In addition, we present the local well-posedness result for the undamped equation in the energy-critical case. For the damped equation, since the coefficient \(\mathrm{e}^{-2\sigma\mathbf{A}(t)}\) after the above transformation is a bounded function, the two equations are almost identical as far as local well-posedness is concerned.
\begin{theorem}[Local well-posedness]\label{局部适定性}
Let \(N\ge 2\), \(\frac{N}{2N-1}<\lambda\le s\le 1\), and \(0<b<\min\{N-2\lambda,2s\}\).
If \(1<2\sigma+1<p_\lambda^c\) (i.e., in the \(H^\lambda\)-subcritical regime) and the initial data \(u_0\in H_{rd}^\lambda\), then there exist a time
\[
T:=T_{N,b,s,\sigma,\|u_0\|_{H^\lambda}}>0
\]
and a unique local solution \(u\) of the inhomogeneous fractional nonlinear Schr\"odinger equation (IFNLS), which belongs to the space
\[
C([0,T],H_{rd}^\lambda)\bigcap_{(q,r)\in\Lambda} L_T^q(W^{\lambda,r}).
\]

Furthermore, if \(2\sigma+1>p_c\) and \(\|u_0\|_{H_{rd}^\lambda}\le A\), then there exists \(\delta:=\delta(A)>0\) such that if
\[
\|e^{-i(-\Delta)^s}u_0\|_{\Omega_{s_c}}\le\delta,
\]
then the solution of problem (IFNLS) is global and satisfies
\begin{equation}
\|u\|_{\Omega_{s_c}}\le 2\delta
\quad\text{and}\quad
\|u\|_{\Omega}+\|(-\Delta)^{\frac{\lambda}{2}}u\|_{\Omega}<\infty.
\end{equation}
Here \(p_\lambda^c:=1+\frac{2(2s-b)}{N-2\lambda}\), \(p_c:=p_0^c=1+\frac{2(2s-b)}{N}\).
\end{theorem}
\begin{remark}
   As the local well-posedness is inherited from the undamped case, the global well-posedness for small initial data in the presence of time-dependent damping follows immediately.
\end{remark}
Below we present the sharp Gagliardo--Nirenberg interpolation inequality with weight \(|x|^{-b}\).
\begin{theorem}\cite{Bhimani2023}
Let \(0<\sigma<\frac{2s-b}{N-2s}\) and \(0<b<\min\{2s,N\}\). Then the sharp Gagliardo--Nirenberg inequality holds:
\[
\int_{\mathbb{R}^N}|x|^{-b}|u|^{2\sigma+2}\,dx
\le K_{\mathrm{opt}}
\|(-\Delta)^{\frac{\sigma}{2}}u\|_{L^2}^{\frac{\sigma N+b}{s}}
\|u\|_{L^2}^{(2\sigma+2)-\frac{\sigma N+b}{s}},
\]
and the sharp constant \(K_{\mathrm{opt}}\) has the explicit expression
\begin{equation}
K_{\mathrm{opt}}
=
\left(
\frac{\sigma N+b}{2s(\sigma+1)-(\sigma N+b)}
\right)^{\frac{2s-(\sigma N+b)}{2s}}
\frac{2s(\sigma+1)}{(N\sigma+b)\|Q\|_{L^2}^{2\sigma}},
\end{equation}
where \(Q\) is the unique nonnegative, radially symmetric, decreasing solution of the equation
\begin{equation}
-(-\Delta)^s Q-Q+|x|^{-b}|Q|^{2\sigma}Q=0.
\end{equation}
Moreover, the solution \(Q\) satisfies the following relations:
\begin{equation}
\|(-\Delta)^{\frac{\sigma}{2}}Q\|_{L^2}
=
\left(
\frac{\sigma N+b}{2s(\sigma+1)-(\sigma N+b)}
\right)^{\frac{1}{2}}
\|Q\|_{L^2},
\end{equation}
and
\begin{equation}
\int_{\mathbb{R}^N}|x|^{-b}|Q|^{2\sigma+2}\,dx
=
\frac{2s(\sigma+1)}{2s(\sigma+1)-(\sigma N+b)}
\|Q\|_{L^2}^2.
\end{equation}
\end{theorem}
Finally, we give an important elementary lemma that can be found in many papers; see, e.g., \cite{Majdoub2023}.
\begin{lemma} \label{一个有用的引理}
Let \(\mathbf{I}\subset\mathbb{R}\) be a time interval, and let \(\mathbf{X}:\mathbf{I}\to[0,\infty)\) be a continuous function satisfying, for every \(t\in\mathbf{I}\),
\[
\mathbf{X}(t)\le a+b[\mathbf{X}(t)]^\theta,
\]
where \(a,b>0\) and \(\theta>1\) are constants. Suppose that for some \(t_0\in\mathbf{I}\),
\[
\mathbf{X}(t_0)\le a,\qquad
a b^{\frac{1}{\theta-1}}
<
(\theta-1)\theta^{\frac{1}{1-\theta}}.
\]
Then, for every \(t\in\mathbf{I}\), we have
\[
\mathbf{X}(t)<\frac{\theta a}{\theta-1}.
\]
\end{lemma}
\begin{lemma}[Fractional Leibniz rule]\label{分数阶1}
Let \(s\in(0,1]\), and let
\[
1<r,r_1,r_2,p_1,p_2<+\infty
\]
satisfy
\[
\frac{1}{r}=\frac{1}{r_i}+\frac{1}{p_i},
\qquad i=1,2.
\]
Then the following estimate holds:
\begin{equation}
\begin{aligned}
\|D^s(fg)\|_{L^r}
\leq {}&
C\|f\|_{L^{r_1}}\|D^s g\|_{L^{p_1}}\\
&+C\|D^s f\|_{L^{r_2}}\|g\|_{L^{p_2}}.
\end{aligned}
\end{equation}
\end{lemma}

\begin{proof}
See the work of Kenig--Ponce--Vega \cite{Kenig1993}.
\end{proof}

\begin{lemma}[Fractional chain rule estimate]\label{分数阶2}
Assume that \(G\in C^1(\mathbb{C})\), \(s\in(0,1]\), and
\[
1<r,r_1,r_2<+\infty
\]
satisfy
\[
\frac{1}{r}=\frac{1}{r_1}+\frac{1}{r_2}.
\]
Then the following estimate holds:
\begin{equation}
\|D^s G(u)\|_{L^r}
\leq C\|G'(u)\|_{L^{r_1}}
\|D^s u\|_{L^{r_2}}.
\end{equation}
\end{lemma}

\begin{proof}
See Kenig--Ponce--Vega \cite{Kenig1993}.
\end{proof}

\section{Proof of the Main Results}
Now we begin the proof of Theorem \ref{main result}.

\begin{proof}
    By the local well-posedness result for the equation, it is enough to show that the \(H^s\) norm of \(u\) is uniformly bounded with respect to \(t\). In the defocusing case the theorem is obvious; hence we only discuss the focusing case below.
    
   Combining the expression for the energy with the sharp Gagliardo--Nirenberg inequality, we obtain
\begin{equation*}
    \begin{aligned}
        E(u(t))
        &=\|(-\Delta)^{\frac{s}{2}}u\|^2_{L^2}
        -\frac{1}{\sigma+1}\int_{\mathbb{R}^N}|x|^{-b}|u(t,x)|^{2\sigma+2}\,dx \\
        &\ge \|(-\Delta)^{\frac{s}{2}}u\|^2_{L^2}
        -\frac{1}{\sigma+1}K_{\mathrm{opt}}
        \|(-\Delta)^{\frac{s}{2}}u\|_{L^2}^{\frac{\sigma N+b}{s}}
        \|u\|_{L^2}^{(2\sigma+2)-\frac{\sigma N+b}{s}}.
    \end{aligned}
\end{equation*}
    Since the mass satisfies \eqref{质量耗散}, we have \(\|u\|_{L^2}\le\|u_0\|_{L^2}\). Substituting this into the above inequality gives
\begin{equation}
    E(u(t))\ge \|(-\Delta)^{\frac{s}{2}}u\|^2_{L^2}
    -\frac{1}{\sigma+1}K_{\mathrm{opt}}
    \|(-\Delta)^{\frac{s}{2}}u\|_{L^2}^{\frac{\sigma N+b}{s}}
    \|u_0\|_{L^2}^{(2\sigma+2)-\frac{\sigma N+b}{s}}.
\end{equation}
Denote
\[
B=\frac{K_{\mathrm{opt}}}{\sigma+1}\|u_0\|_{L^2}^{2\sigma+2-\frac{\sigma N+b}{s}}.
\]
We first show that \(I(u(t))>0\) for every \(t\in\mathbb{R}^+\). Since we assume \(I(u_0)>0\), suppose that
\[
t^*=\inf\{t\in\mathbb{R}^+\mid I(u(t))<0\}.
\]
By the continuity of \(I(u)\) in time, we know that \(t^*>0\) and \(I(u(t^*))=0\). That is, at time \(t^*\), the kinetic and potential energies must be equal.

Using the definition of \(t^*\), we know that \(I(u(t))>0\) for every \(t\in[0,t^*)\). According to \eqref{能量导数}, \(E(u)\) is monotonically decreasing on \([0,t^*)\). Hence \(E(u(t^*))\le E(u_0)<E(Q)\). However, when the kinetic energy equals the potential energy, the energy of the solution cannot be lower than the ground-state energy; this follows from the definition of the ground state and the variational principle. Namely, \(E(u(t^*))\ge E(Q)\), a contradiction. Therefore \(t^*=\infty\). Consequently, for every \(t\in\mathbb{R}^+\), \(E(u(t))\) is monotonically decreasing and \(E(u(t))\le E(u_0)\). Combining this with the above inequality, we have
   \begin{equation}
    E(u_0)\ge \|(-\Delta)^{\frac{s}{2}}u\|^2_{L^2}
    -\frac{1}{\sigma+1}K_{\mathrm{opt}}
    \|(-\Delta)^{\frac{s}{2}}u\|_{L^2}^{\frac{\sigma N+b}{s}}
    \|u_0\|_{L^2}^{(2\sigma+2)-\frac{\sigma N+b}{s}}
    =\|(-\Delta)^{\frac{s}{2}}u\|^2_{L^2}
    -B\|(-\Delta)^{\frac{s}{2}}u\|^{\frac{\sigma N+b}{s}}_{L^2}.
\end{equation}
Let \(f(X(t))=X(t)-B X(t)^\theta\), where \(\theta=\frac{\sigma N+b}{2s}\). It is easy to verify that the conditions of Lemma \ref{一个有用的引理} are satisfied, and hence
\[
X(t)<\frac{\theta E(u_0)}{\theta-1}=C.
\]
Therefore, \(\|u\|_{H^s}\) is uniformly bounded for all \(t\in\mathbb{R}^+\), and thus \(u\) exists globally.

Next, we give an important estimate of the nonlinear term.
\begin{lemma} \label{非线性估计}
   Let \(N\ge 3\). Under the conditions on the parameters in Theorem \ref{局部适定性}, for the nonlinear term we have the following estimate:
\begin{equation}
\left\|\nabla^s\left(e^{-2\sigma A(t)}|x|^{-b}|u|^{2\sigma}u\right)\right\|_{\Omega'(I)}
\le
Ce^{-C_1t_0}
\|u\|_{L_t^\infty H_x^s}^{\theta}
\|\nabla^s u\|_{\Omega(I)}^{2\sigma+1-\theta},
\end{equation}
where \(\theta\) is a sufficiently small positive constant, and the estimate holds under the large damping condition, namely \(\underline{a}>0\); the interval \(I\) denotes \([t_0,\infty)\).
\end{lemma}
\begin{proof}
    We take
\[
r=\frac{(4\sigma+4-2\theta)N(2N-1)}{2N(2N-1)+4\sigma N-(\theta N+2b)(2N-1)},
\qquad
q=\frac{2s(2\sigma+2-\theta)(2N-1)}{2N^2-3\sigma N+2bN-b}.
\]
It is easy to verify that \((q,r)\in\Lambda\).

First, we estimate
\[
\left\|\nabla^s\left(e^{-2\sigma A(t)}|x|^{-b}|u|^{2\sigma}u\right)\right\|_{L^{r'}}.
\]
We denote the unit ball \(B_1(0)=B\). Then
\begin{equation}
\left\|\nabla^s\left(e^{-2\sigma A(t)}|x|^{-b}|u|^{2\sigma}u\right)\right\|_{L^{r'}}
=
e^{-2\sigma A(t)}
\left\|\nabla^s\left(|x|^{-b}|u|^{2\sigma}u\right)\right\|_{L^{r'}}
\le
e^{-2\sigma A(t)}\bigl(N(t,B)+N(t,B^C)\bigr),
\end{equation}
where
\[
N(t,B)=\left\|\nabla^s\left(|x|^{-b}|u|^{2\sigma}u\right)\right\|_{L^{r'}(B)},
\qquad
N(t,B^C)=\left\|\nabla^s\left(|x|^{-b}|u|^{2\sigma}u\right)\right\|_{L^{r'}(B^C)}.
\]
By the fractional differentiation rules (Lemmas \ref{分数阶1}, \ref{分数阶2}), we have
\[
N(t,B)
=
\left\||x|^{-b-s}|u|^{2\sigma}u\right\|_{L^{r'}(B)}
+
\left\||x|^{-b}|u|^{2\sigma}\nabla^s u\right\|_{L^{r'}(B)}
=
N_1(t,B)+N_2(t,B).
\]
We first consider \(N_1(t,B)\). Using the Sobolev embedding
\begin{equation}
\begin{aligned}
    N_1(t,B) &= \left\| |x|^{-b-s} |u|^{2\sigma} u \right\|_{L^{r'}(B)} \\
    &\le \left\| |x|^{-b-s} \right\|_{L^\gamma(B)} \left\| u^\theta \right\|_{L^a} \left\| u^{2\sigma+1-\theta} \right\|_{L^c} \\
    &= \left\| |x|^{-b-s} \right\|_{L^\gamma(B)} \left\| u \right\|_{L^{a\theta}}^\theta \left\| u \right\|_{L^{c(2\sigma+1-\theta)}}^{2\sigma+1-\theta} \\
    &\le \left\| |x|^{-b-s} \right\|_{L^\gamma(B)} \left\| u \right\|_{L^{a\theta}}^\theta \left\| \nabla^s u \right\|_{L^r}^{2\sigma+1-\theta},
\end{aligned}
\end{equation}
where the exponents satisfy
\begin{equation*} 
    \begin{aligned}
        & a\theta=\frac{2N}{N-2s},\\
        & \frac{1}{r}-\frac{1}{c(2\sigma+1-\theta)}=\frac{s}{N},\\
        & \frac{1}{r'}=\frac{1}{\gamma}+\frac{1}{a}+\frac{1}{c}.
    \end{aligned}
\end{equation*} 
Hence we have
\begin{equation}
    \frac{N}{\gamma}=N-\frac{N}{r}-\frac{\theta(N-2s)}{2}-\left(\frac{N}{r}-s\right)(2\sigma+1-\theta).
\end{equation}
Substituting the value of \(r\), we obtain
\begin{equation*}
    \frac{N}{\gamma}=b+s+2\sigma\left(s-\frac{N}{2N-1}\right)>b+s.
\end{equation*}
Therefore,
\begin{equation}
    \left\||x|^{-b-s}\right\|_{L^{\gamma}(B)}\le C.
\end{equation}
Combining the above estimates, we get
\begin{equation}
    N_1(t,B)\le C\|u\|^\theta_{H^s}\|\nabla^s u\|^{2\sigma+1-\theta}_{L^r}.
\end{equation}
Next, we consider \(N_2(t,B)\). Similarly, using Hölder's inequality together with the Sobolev embedding, we obtain
\begin{equation}
\begin{aligned}
    N_2(t,B) &\le \left\| |x|^{-b} \right\|_{L^{\gamma}(B)} \left\| u^\theta \right\|_{L^a} \left\| u^{2\sigma-\theta} \right\|_{L^c} \left\| \nabla^s u \right\|_{L^r} \\
    &= \left\| |x|^{-b} \right\|_{L^{\gamma}(B)} \left\| u \right\|_{L^{a\theta}}^\theta \left\| u \right\|_{L^{c(2\sigma-\theta)}}^{2\sigma-\theta} \left\| \nabla^s u \right\|_{L^r}.
\end{aligned}
\end{equation}
where the exponents satisfy
\begin{equation*}
    \begin{aligned}
        & a\theta=\frac{2N}{N-2s},\\
        & \frac{1}{r}-\frac{1}{c(2\sigma-\theta)}=\frac{s}{N},\\
        & \frac{1}{\gamma}+\frac{1}{a}+\frac{1}{c}+\frac{1}{r}=\frac{1}{r'}.
    \end{aligned}
\end{equation*}
Similarly, it follows from the above that
\begin{equation}
\begin{aligned}
    \frac{N}{\gamma}
    &= N-\frac{2N}{r}-\frac{N}{a}-\left(\frac{N}{r}-s\right)(2\sigma-\theta) \\
    &= N-\frac{\theta N}{2}+2s\sigma-\frac{N}{r}(2\sigma+2-\theta).
\end{aligned}
\end{equation}
Again substituting the value of \(r\), we obtain
\begin{equation}
    \frac{N}{\gamma}=b+2\sigma\left(s-\frac{N}{2N-1}\right)>b.
\end{equation}
Therefore,
\begin{equation}
    \left\||x|^{-b}\right\|_{L^{\gamma}(B)}\le C,
\end{equation}
and hence
\begin{equation}
    N_2(t,B)\le C\|u\|^\theta_{H^s}\|\nabla^s u\|^{2\sigma+1-\theta}_{L^r}.
\end{equation}
Next, we consider \(N(t,B^c)\). Similarly to the preceding argument, we decompose it as follows:
\begin{equation}
    N(t,B^c)=N_3(t,B^c)+N_4(t,B^c)
    =\left\||x|^{-b-s}|u|^{2\sigma}u\right\|_{L^{r'}(B^c)}
    +\left\||x|^{-b}|u|^{2\sigma}\nabla^s u\right\|_{L^{r'}(B^c)}.
\end{equation}
For \(N_3(t,B^c)\), we estimate as follows:
\begin{equation}
\begin{aligned}
    N_3(t,B^c) &= \left\| |x|^{-b-s} |u|^{2\sigma} u \right\|_{L^{r'}(B^c)} \\
    &\le \left\| |x|^{-b-s} \right\|_{L^\gamma(B^c)} \left\| u^\theta \right\|_{L^a} \left\| u^{2\sigma+1-\theta} \right\|_{L^c} \\
    &= \left\| |x|^{-b-s} \right\|_{L^\gamma(B^c)} \left\| u \right\|_{L^{a\theta}}^\theta \left\| u \right\|_{L^{c(2\sigma+1-\theta)}}^{2\sigma+1-\theta} \\
    &\le \left\| |x|^{-b-s} \right\|_{L^\gamma(B^c)} \left\| u \right\|_{L^{a\theta}}^\theta \left\| \nabla^s u \right\|_{L^r}^{2\sigma+1-\theta}.
\end{aligned}
\end{equation}
The exponents satisfy the following conditions:
From these relations, we have
\begin{equation}
    \frac{N}{\gamma}-b-s=-\frac{N}{a}+2\sigma\left(s-\frac{N}{2N-1}\right)+\theta\left(\frac{N}{2}-s\right).
\end{equation}
Since \(\frac{N}{2N-1}<s\le1<\frac{N}{2}\), we set \(2\sigma(s-\frac{N}{2N-1})=d>0\). Therefore, as long as we take
\[
a\theta<\frac{2N}{\frac{2d}{\theta}+N-2s}
\quad\text{and}\quad
\theta\ge\frac{d}{s},
\]
we have \(\frac{N}{\gamma}-b-s<0\), and
\[
2\le a\theta<\frac{2N}{\frac{2d}{\theta}+N-2s}<\frac{2N}{N-2s}.
\]
Hence, in this case we obtain
\begin{equation}
    N_3(t,B^c)\le C\|u\|^\theta_{H^s}\|\nabla^s u\|_{L^r}^{2\sigma+1-\theta}.
\end{equation}
Finally, we turn to \(N_4(t,B^c)\). Similarly, we have
\begin{equation}
\begin{aligned}
    N_4(t,B^c) &= \left\| |x|^{-b} |u|^{2\sigma} \nabla^s u \right\|_{L^{r'}(B^c)} \\
    &\le \left\| |x|^{-b} \right\|_{L^{\gamma}(B^c)} \left\| |u|^{\theta} \right\|_{L^a} \left\| |u|^{2\sigma-\theta} \right\|_{L^c} \left\| \nabla^s u \right\|_{L^r} \\
    &= \left\| |x|^{-b} \right\|_{L^\gamma(B^c)} \left\| u \right\|_{L^{a\theta}}^\theta \left\| u \right\|_{L^{c(2\sigma-\theta)}}^{2\sigma-\theta} \left\| \nabla^s u \right\|_{L^r} \\
    &\le \left\| |x|^{-b} \right\|_{L^\gamma(B^c)} \left\| u \right\|_{L^{a\theta}}^\theta \left\| \nabla^s u \right\|_{L^r}^{2\sigma+1-\theta}.
\end{aligned}
\end{equation}
The exponents satisfy
\begin{equation}
    \begin{aligned}
        &\frac{1}{\gamma}+\frac{1}{a}+\frac{1}{c}+\frac{1}{r}=\frac{1}{r'},\\
        &\frac{1}{r}-\frac{1}{c(2\sigma-\theta)}=\frac{s}{N},\\
        & 2\le a\theta<\frac{2N}{N-2s}.
    \end{aligned}
\end{equation}
Therefore, we can obtain
\begin{equation}
    \frac{N}{\gamma}-b
    =-\frac{N}{a}+2\sigma\left(s-\frac{N}{2N-1}\right)+\theta\left(\frac{N}{2}-s\right)
    =-\frac{N}{a}+d+\theta\left(\frac{N}{2}-s\right).
\end{equation}
Similarly, choosing
\[
a\theta<\frac{2N}{\frac{2d}{\theta}+N-2s},
\qquad
\theta\ge\frac{d}{s},
\]
we obtain
\begin{equation}
    N_4(t,B^c)\le C\|u\|^\theta_{H^s}\|\nabla^s u\|_{L^r}^{2\sigma+1-\theta}.
\end{equation}
In summary,
\begin{equation}
\begin{aligned}
    \left\| \nabla^s \left( e^{-2\sigma A(t)} |x|^{-b} |u|^{2\sigma} u \right) \right\|_{L^{r'}} 
    &= e^{-2\sigma A(t)} \left\| \nabla^s \left( |x|^{-b} |u|^{2\sigma} u \right) \right\|_{L^{r'}} \\
    &\le e^{-2\sigma A(t)} \left( N(t,B) + N(t,B^c) \right) \\
    &\le C e^{-2\sigma A(t)} \left\| u \right\|_{H^s}^\theta \left\| \nabla^s u \right\|_{L^r}^{2\sigma+1-\theta}.
\end{aligned}
\end{equation}
The next task is to take the norm in time:
\begin{equation}
\begin{aligned}
    &\left\| \left\| \nabla^s \left( e^{-2\sigma A(t)} |x|^{-b} |u|^{2\sigma} u \right) \right\|_{L^{r'}} \right\|_{L^{q'}} \\
    &\le C \left\| e^{-2\sigma A(t)} \right\|_{L^{\frac{q}{q-2\sigma-2+\theta}}} \left\| u \right\|_{L^\infty(H^s)}^\theta \left\| \left\| \nabla^s u \right\|_{L^r}^{2\sigma+1-\theta} \right\|_{L^{\frac{q}{2\sigma+1-\theta}}} \\
    &\le C \left\| e^{-2\sigma A(t)} \right\|_{L^{\frac{q}{q-2\sigma-2+\theta}}} \left\| u \right\|_{L^\infty(H^s)}^\theta \left\| \nabla^s u \right\|_{\Omega}^{2\sigma+1-\theta}.
\end{aligned}
\end{equation}
Substituting the value of \(q\), one finds that \(\frac{q}{q-2\sigma-2+\theta}\) is greater than or equal to \(1\); we denote it by \(\tau\). Therefore,
\begin{equation}
    \|e^{-2\sigma A(t)}\|^{\tau}_{L^\tau(I)}
    \le
    \int_{t_0}^{\infty} e^{-c\underline{\mathbf{a}}t\tau}\,dt
    =
    C_1e^{-C_2t_0}
\end{equation}
The constants appearing before and after may denote different quantities; the lemma is proved.
\end{proof}

By the preceding nonlinear estimate lemma \ref{非线性估计}, combined with the Strichartz estimates, we obtain
    \begin{equation}
\begin{aligned}
    \left\| \nabla^s u \right\|_{L^q(I, L^r)} 
    &\lesssim \left\| u_0 \right\|_{H^s} + \left\| \nabla^s \left[ e^{-2\sigma A(t)} |x|^{-b} |u|^{2\sigma} u \right] \right\|_{L^{q'}(I, L^{r'})} \\
    &\lesssim \left\| u_0 \right\|_{H^s} + C e^{-C_1 t_0} \left\| u \right\|_{L_t^\infty H_x^s}^\theta \left\| \nabla^s u \right\|_{\Omega(I)}^{2\sigma+1-\theta}.
\end{aligned}
\end{equation}
Here the interval is \(I=[t_0,\infty)\), and \((q,r)\) are the exponents from Lemma \ref{非线性估计}. As shown above, \(\|u\|_{H^s}\) is uniformly bounded in time; hence we can choose \(t_0\) sufficiently large so that \(Ce^{-C_1t_0}\|u\|_{L^{\infty}(H^s)}\) is sufficiently small. Then, applying Lemma \ref{一个有用的引理} again, or using a standard continuity argument, we obtain
\begin{equation}
    \|u\|_{L^q(I,L^r)}\lesssim \|u_0\|_{H^s}.
\end{equation}
Next, we turn to scattering. Given \(t_1>t_2>t_0\), by the Strichartz estimates we obtain
   \begin{equation}
\begin{aligned}
    \left\| e^{it_1(-\Delta)^s} u(t_1) - e^{it_2(-\Delta)^s} u(t_2) \right\|_{H^s}
    &\le \left\| \int_{t_2}^{t_1} e^{i\tau(-\Delta)^s} e^{-2\sigma A(\tau)} \left[ |x|^{-b} |u|^{2\sigma} u \right] \mathrm{d}\tau \right\|_{H^s} \\
    &\lesssim \left\| \nabla^s u \right\|_{\Omega(t_1,t_2)}^{2\sigma+1-\theta} \rightarrow 0
\end{aligned}
\end{equation}
as \(t_1,t_2\rightarrow\infty\). Hence \(e^{it(-\Delta)^s}u(t)\) is a Cauchy sequence in the space \(H^s\). By the completeness of \(H^s\), we obtain the limit
\begin{equation}
    u_+ := u_0 - i\mu \int_{0}^{\infty} e^{is(-\Delta)^s} e^{-2\sigma A(s)} |v(s)|^{2\sigma} v(s)\,ds.
\end{equation}
It exists in \(H^s\). Repeating the above argument, we can prove that
\begin{equation}
\lim_{t \to \infty} \|u(t) - e^{-it(-\Delta)^s}u_+\|_{H^1} = 0,
\end{equation}
thus completing the proof of the theorem.
\end{proof}

\end{document}